\documentclass[10pt]{article}
\usepackage{amsmath,amssymb,amsthm,amscd}
\numberwithin{equation}{section}

\newtheorem{prop}{Proposition}[section]
\newtheorem{theo}[prop]{Theorem}
\newtheorem{lem}[prop]{Lemma}
\newtheorem{cor}[prop]{Corollary}
\newtheorem{rem}[prop]{Remark}

\newtheorem{defi}[prop]{Definition}
\newtheorem{conj}[prop]{Conjecture}
\newtheorem{q}[prop]{Question}

\def\begeq{\begin{equation}}
\def\endeq{\end{equation}}

\def\and{\quad{\rm and}\quad}

\let\lra=\longrightarrow

\def\mapright\#1{\,\smash{\mathop{\lra}\limits^{\#1}}\,}

\begin{document}

\title{Calabi flow in Riemann surfaces revisited: A new point of view}
\author{ X. X. Chen\footnote{Research was supported partially by NSF postdoctoral
fellowship.}}
\date{ March 12, 1999}
\pagestyle{plain}
\bibliographystyle{plain}

\maketitle
\section{Introduction}
   P. Chrusciel's now famous paper \cite{Chru91} shows that the Calabi flow exists
and converges to a constant scalar curvature metric
  in a Riemann surface,
assuming the existence of a constant scalar curvature metric in
the background. Since there always exists a constant scalar
curvature metric
 thanks to the uniformization theorem in
Riemannian surface, Chrusciel's proof appears to be satisfactory for
most purposes.  One memorable feature
of Chrusciel's paper is the strong and somewhat mysterious
influence of physics. Given the importance of the Calabi flow in K\"ahler geometry,
  a more direct mathematical proof of Chrusciel's theorem is desirable. Inspired by Chrusciel's work, the author
has worked on a related problem on the existence of extremal
K\"ahler metrics  in any
 Riemannian surface with boundary since 1994 (cf.  \cite{chenthesis94}, \cite{chen942}, \cite{chen981} for further references).
The local theory we developed in the aforementioned
 papers together with our observation that the Calabi flow decreases several interesting  functionals (cf. Section 3),
 are sufficient to provide a new proof of Chrusciel's theorem.
In this paper, we also need to assume the uniformization theorem.
At the end
of this paper, we will remark  on how to remove this assumption.\\

  The purpose of this short note is two-fold: first, in
Chrusciel's original paper, the Bondi mass estimate plays a
crucial role. The idea of Bondi mass might be  clear to  people
with some physics background, but it is hard to connect to those
of us  who are less well versed in physics. It is also difficult
to find a higher dimensional analogue of ``Bondi mass.''
 Thus,  a proof without the
Bondi mass estimate would be appreciated by mathematicians. The
second, and the most important
 point is to cast this Calabi flow
problem into new light. In numerous instances in the history of
mathematics,
 a new proof to an old problem from a completely
different angle  has lead to progress in other, often unrelated problems.
 In general, for a parabolic equation, one first proves the
long time existence of the flow; then one argues that the flow
converges by sequence of times $t_i \rightarrow \infty$. The
uniqueness of the limit for different sequences as time
approaches $\infty$ is  usually one of the hardest problems.
Moreover, there are very few tools available (cf.  \cite{Lms83})
to tackle this problem. The recent development on a Riemannian
metric in  the space of K\"{a}hler metrics (cf.  \cite{Ma87},
\cite{Semmes92} and \cite{Dona96} for more references) provides
more tools for us to argue on this point\footnote{For convenience
of the readers, we briefly discuss  these new developments in
section 2.}. Let $\cal H$ denote the space of K\"{a}hler
potentials in a fixed K\"ahler class (any K\"ahler metric in a
fixed K\"ahler class determines a unique K\"ahler
 potential
up to some additive constant. Thus, we sometimes use the term
 "space of K\"ahler metrics" and  "space of K\"ahler potentials"
interchangeably.) Following a program outlined by  Donaldson
\cite{Dona96}, we proved in \cite{chen981}  that the $\cal H$
 is geodesically convex by $C^{1,1}$ geodesic and more importantly it is a
metric space (cf. Theorem B and Definition 2 in Section 2).
\begin{defi} (Cauchy curve): Any curve $c(t), 0\leq t < \infty$ is
called a Cauchy curve in $\cal H$ if for any $t,s \rightarrow \infty$
the  distance of $c(t)$ and $c(s)$ in $\cal H$ approaches
$0$ uniformly.
\end{defi}
  For any initial metric, we prove that the curve in $\cal H$
obtained by the Calabi flow is precisely a Cauchy curve in ${\cal
H}.\;$ The uniqueness of the limit by sequences is an easy
consequence of
this assertion. \\

{\bf Organization:}  In Section 3, we describe four energy
functionals, each of which is non-increasing under the Calabi
flow. Lower bounds for these functionals give us some useful
integral estimates.  In Section 4, we prove the long time
existence of the Calabi flow on a Riemann surface via these
integral estimates. In Section 5, we prove that for any sequence
of $t_i \rightarrow \infty,$ there exists a subsequence where the
Calabi flow converges to a constant scalar curvature metric. In
Section 6, we use a new technique introduced in Section 2 to prove
the uniqueness of the sequential limit of the Calabi flow as $t
\rightarrow \infty.$ Moreover, we also give an analytic proof of
convergence of the
 Calabi flow  to
a constant scalar curvature metric (We don't need to adjust the
flow by some conformal transformation).
 In Section 7,
we discuss some interesting questions for future study.\\

  The bulk of this work was carried out when the author was
visiting Stanford University in 1998. Thanks to Stanford
university for their hospitality,   and   Rick Schoen and Leon
Simon for many interesting discussions. Thanks also to  Professor
Donaldson; his encouragement and interest is the key to this new
proof. The author also wants to thank Professor Xu Xingwang for
pointing out some errors in an earlier version of this paper, and
thank Professor M. Struwe for pointing a computational error in
Section 6 (for proving the exponential convergence)
 in the earlier version. Thanks also to Wang Guofang for carefully reading an earlier version of this paper
 and for many interesting suggestions.  Finally, I am particularly grateful for the referee's insightful comments.\\

\section{Summary of recent developments
in the space of K\"ahler potentials} Let $(V,\omega_0)$ be an
n-dimensional K\"ahler manifold.
 Mabuchi (\cite{Ma87})  in 1987 defined a Riemannian metric
on the space of K\"ahler metrics,
under which it  becomes (formally) a non-positive curved infinite dimensional
 symmetric space. Apparently unaware of Mabuchi's work,
Semmes \cite{Semmes92}  and Donaldson \cite{Dona96} re-discovered
this same metric again from different angles. Set
\[
 {\cal H} = \{ \varphi | \omega_{\varphi} = \omega_0 + \sqrt{-1} \partial \overline{\partial} \varphi > 0 \;\;{\rm on} \;V\}.
\]
Clearly, the tangent space of $\cal H$ is $C^{\infty}(V)$ if we assume that everything is smooth.
For any  vector $\psi$ in the tangent space $ T_{\varphi} \cal {H}, $ we define the length of this vector as
\[
\|\psi\|^2_{\varphi} =\int_{V}\psi^2\;{\omega_{\varphi}}^n \;
d\,\mu_{\varphi},
\]
where $d\;\mu_{\varphi}$ is the volume element of the corresponding K\"ahler metric determined by $\varphi.$
The geodesic equation is
\begin{equation}
  \varphi(t)'' - {1\over 2} |\nabla \varphi'(t)|^2_{ \varphi(t)} = 0,
\label{geodesic}
\end{equation}
where the derivative and norm in the second term of the left hand
side
are taken with respect to the metric $\omega_{\varphi(t)}.\;$\\

This geodesic equation shows us how to define a connection on the
tangent bundle of ${\cal H}$.  If $\phi(t)$ is any path in ${\cal
H}$ and $\psi(t)$ is a field of tangent vectors along the path
(that is, a function on $V \times [0,1]$), we define the
covariant derivative along the path to be
$$ D_{t}\psi = \frac{\partial\psi}{\partial t} - {1\over 2}  (\nabla \psi, \nabla \phi')_{\phi}.  $$
 The main theorem formally
proved  in \cite{Ma87} (and later re-proved in \cite{Semmes92} and
 \cite{Dona96}) is:\\

\noindent {\bf Theorem A} {\it The Riemannian manifold $\cal {H} $ is an infinite dimensional symmetric space; it admits a Levi-Civita connection whose
curvature is covariant constant. At a point $\phi\in{\cal {H}}$ the curvature  is given by
\[     R_{\phi}(\delta_{1}\phi, \delta_{2}\phi) \delta_{3}\phi=
- {1\over 4}  \{ \{ \delta_{1}\phi, \delta_{2}\phi\}_{\phi},
\delta_{3}\phi\}_{\phi},\]
where $\{\ ,\ \}_{\phi}$ is the Poisson bracket on $C^{\infty}(V)$ of the
symplectic form $\omega_{\phi}$; and $\delta_1 \phi, \delta_2 \phi \in T_{\phi} {\cal H}.\;$ Then the sectional curvature is non-positive, given by}
\[    K_{\phi}(\delta_{1}\phi, \delta_{2}\phi) = - {1\over 4}  \Vert
 \{ \delta_{1} \phi , \delta_{2}\phi\}_{\phi}\Vert_{\phi}^2. \]
We will skip the proof later, interested readers are referred to
\cite{Ma87} or \cite{Semmes92} and
 \cite{Dona96} for the proof.\\

The subject has been quiet since the early pioneer work of Mabuchi
(1987) and Semmes (1991). The real breakthrough came in the
beautiful paper by Donaldson \cite{Dona96} in 1996, where
 he  outlines the connection between this
Riemannian metric in the infinite dimensional space $\cal H$ and
the traditional K\"ahler geometry,  via a series of important
conjectures and theorems. In 1997, following his program, the
author proves some of
his conjectures:\\

\noindent {\bf Theorem B} \cite{chen991}{\it The following statements
are true:
\begin{enumerate}
\item The space of K\"ahler potentials ${\cal H}$ is convex by $C^{1,1}$ geodesic.
More precisely, if $\varphi_0,\varphi_1 \in \cal H,\;$ then there
exists a unique geodesic path $\varphi(t) \;(0\leq t \leq 1)$
connecting these two points, such that the mixed covariant
derivatives of $\varphi(t)$ are uniformly bounded from above.
\item  ${\cal H}$ is a metric space. In other words, the infimum of the
lengths of all possible curves between any two different points in
$\cal H$ is strictly positive.
\end{enumerate}
}

\begin{defi}  For any two K\"{a}hler metrics $g_1$ and $g_2$ in $\cal H,$ define
the distance $d(g_1,g_2)$ to be  the length of geodesic
connecting them.
\end{defi}
\begin{rem} One can complete $\cal H$ by adding  all of the limits
 of the Cauchy sequence under this distance function.
\end{rem}

In \cite{chen992}, E. Calabi and the author proved the following:\\

\noindent {\bf Theorem C}\cite{chen992}{\it The following statements are true:
\begin{enumerate}
\item  $\cal H$ is a non-positive curved space in the sense of Alenxandrov.
\item  The length of any smooth curve in $\cal H$ is decreasing under
the Calabi flow unless it is represented by a holomorphic
transformation. The distance in $\cal H$ is also decreasing if
the Calabi flow exists for all the time (from $t=0$ to $\infty$)
for any initial smooth metric.
\end{enumerate}}

\section {Various energy functionals in Riemann surfaces}
\subsection{Notations}
In this subsection, we set up some notations for later use.
Suppose that $g_0 = \sqrt{-1} F_0 d\,z d\,\overline { z} $ is a
fixed metric and $g = \sqrt{-1} F d\,z d\,\overline { z} $ is any
metric in the same K\"{a}hler class with $g_0.\;$ Then $g$ is
necessarily conformal to $g_0$ with same area.  Suppose now that
$\varphi$ is  K\"{a}hler distortion potential of $g$ w.r.t.
$g_0$, while $e^{2 u}$ is the conformal factor of $g$ w.r.t.
$g_0.$ Then
\[ \sqrt{-1} F  d\,z \wedge d\,\overline { z}=  \sqrt{-1}F_0  d\,z \wedge d\,\overline { z}
 + \sqrt{-1}\partial \overline {\partial} \varphi,\]
and
\[
    g = e^{2 u} g_0 = (1 + \triangle_0 \varphi) g_0.
\]
Here we use $\triangle, \triangle_0, \triangle_g $ to denote the
complex Laplacian operator w.r.t. the local metric $d\,z \wedge
d\,\overline { z}, \; g_0 $ and $g$ respectively. Then \[
\triangle = {{\partial^2}\over {\partial z \overline{\partial
z}}}, \qquad \triangle_g = {\triangle_0 \over e^{2 u}}.
\]

Let $K, K_0$ denote the scalar curvature of $g,\;g_0$
respectively, and let $\underline{K}$ be the average of the scalar
curvature\footnote{We shall normalize so that $\underline{K} $ is
$ 1, 0, -1 $ for $S^2$, torus  and high genus surface.}. Then
\begin{eqnarray}
  K & = & - {{ \triangle (\ln F) }\over {F}} \label{eq:curvature0}\\
    &  = &  {{ -\triangle_0 (\ln (1 + \triangle_0 \varphi)) + K_0 } \over {(1 + \triangle_0 \varphi)}}
     ={{ - 2  \triangle u + K_0}\over {e^{2 u}}}.
\label{eq:curvature}
\end{eqnarray}
The Calabi flow is then defined as
\[
   {{\partial \varphi} \over {\partial t}} = K - \underline{K},
\]
or equivalently
\[
    {{\partial u} \over {\partial t}} = {1\over 2} {{\triangle_0 K}\over{e^{2 u}}} = {1\over 2}\;\triangle_g K.
\]
Let $H^{p,q}(M,g)$ denote the Sobolev space of functions whose
$p-$th derivatives are $L^q$ integral with respect to the metric
$g.\;$ \\

Finally, if $M=S^2, $ let $G$ denote the group of conformal
transformations. Then $ {\rm dim}\; G = 3.\;$ Use $\eta(S^2)$
denote its Lie algebra, or the space of all holomorphic vector
fields in $S^2.\;$

\subsection {Various energy functionals}
  There are several energy functionals which are decreasing under
   the Calabi flow. This is quite unusual since most heat flows
will just decrease one energy functional. In this subsection, we
will  introduce four of those functionals (this is not a complete
list).

\begin{enumerate}
\item {\bf Area functional}:$\;$ Area functional is
 invariant under the Calabi flow.
Denote area functional $A(g)$ as $A = \int_M d\,g.\;$ Then
\[
 {{d\, A(g)}\over {d\,t}} =2 \int_M {{\partial u}\over {\partial t}} d\,g = \int_M \triangle_g K d\,g = 0. \]

\item {\bf The Calabi energy} was first introduced
by E. Calabi in his famous paper \cite{calabi82}. This functional
 measures essentially the $L^2$ distance of a given metric from being a constant scalar curvature metric. Namely,
  \[ Ca(g) = \int_M \; (K - \underline{K})^2 d\, g
\]
or
\begin{equation}
  {d\over {d\,t}} Ca(g(t)) = - \int_M \overline{L}^* L (K) {{\partial \varphi}\over {\partial t}} d\,g,  \label{eq:Calabi-energy}
\end{equation}
where $L$ is the Lichernowicz operator w.r.t. the metric $g,\;
i.e.,$ for any function $f$ on $M,$
\[  L (f) = f_{,zz} d\,z\otimes d\,z.
\]
Here $f_{,zz}$ is the second covariant derivatives of $f$ of
 pure type:
 \begin{equation}
    f_{,zz} = {{\partial^2 f}\over {\partial z^2}} - {{\partial f}\over {\partial
    z}} {{\partial \log F}\over {\partial z}}.
    \label{eq:fzzdefinition}
 \end{equation}
And $\overline{L}^* L$ is a fourth order self adjoint operator
with nonnegative eigenvalue. An important problem in K\"{a}hler
geometry is to estimate the lower bound of the first non-zero
eigenvalue of this 4th order operator over a family of K\"{a}hler
metrics (such as the family produced by Calabi flow). Under the
Calabi flow, we have
\[
\begin{array}{lcl}
  {d\over {d\,t}} Ca(g) & = & - \int_M \overline{L}^* L (K) \;(K - \underline{K}) d\,g \\
   & = & - \int_M |L (K)|^2_g\; d\,g \leq 0.
   \end{array}\]
The last equality holds unless $L (K) \equiv 0, $ or $ g$ is an extremal
K\"{a}hler metric in Riemann surface. By a theorem of Calabi,  $K$ must
be a constant everywhere in this case. Therefore, during the flow, the Calabi
energy must be strictly decreasing.

\item {\bf The Mabuchi energy} was first introduced by Mabuchi in 1987
to prove the uniqueness of K\"{a}hler-Einstein metric when the
first Chern class is positive \cite{Bando87}. It is also
``famous''
 since it is only defined through its derivatives and it is hard
to  manipulate directly. For any smooth curve $\varphi(t) \in
\cal H,\;$ the Mabuchi energy is defined by
\begin{equation}
{d\over {d\,t}} Ma(g) = - \int_M (K - \underline{K}) {{\partial \varphi}\over {\partial t}} d\,g.  \label{eq:mabuchi-energy}
\end{equation}
Note that the Calabi flow is the $L^2$ gradient flow of the
Mabuchi energy.  Consequently, it is also decreasing under the
Calabi flow:
\begin{equation}
  {d\over {d\,t}} Ma(g) = - \int_M (K - \underline{K})^2 d\,g.
  \label{eq:mabuchi-energy1}
\end{equation}

  In a Riemann surface, since there always exists a constant scalar curvature metric, the Mabuchi energy
  has a uniform lower bound (c.f., \cite{Bando87}).
  Then (denote the metric along Calabi flow as $g(t)$):
\[
   Ma(g(t)) |_0^{T} = - \int_0^T \int_M (K - \underline{K})^2 d\,g d\,t > - C.
\]
The constant $C > 0 $ is independent of the time $T.\;$ In other words,
if the flow exists for long time, we have:
\begin{equation}
  \int_0^{\infty} Ca(g(t)) d\,t = \int_0^{\infty} \int_M (K - \underline{K})^2 d\,g d\,t < C.
\end{equation}
 Since we also know that  the Calabi energy itself is decreasing,
this  implies that Calabi energy decreases to $0$ as $t$
approaches $\infty.\;$

\item {\bf The Liouville energy} \footnote{According to \cite{OPS88}, the Liouville
energy represents the log determinant of the Laplacian operator
of conformal metrics $g.\;$ The research in this direction is
very active, see \cite{OPS88}, \cite{ChangQ96}, \cite{ChangY95}
for further references in this topic.}. From the definition of
the Calabi energy and the Mabuchi energy, especially the
definitions via their first derivatives, I believe that there
should exist a third functional $F$ whose definition involves
some second order operator $\cal P$ such that
\[
  {d\over {d\,t}} F(g) = - \int_M  {\cal P} (K) {{\partial \varphi}\over {\partial t}} d\,g,
\]
where
\[
  {\cal P} = \triangle_g +  {\rm lower\; order\; terms}
\]
and $\triangle_g$ is the complex Laplacian operator w.r.t. metric
$g.\;$ In dimension 2, it turns out  that $ {\cal P} =
\triangle_g $ is precisely what it needs, and the functional $F$
above  is just the well-known Liouville functional. It will be
interesting if one could generalize this to higher dimensional
K\"{a}hler manifolds.

\[ F  = \int_M \ln {g \over {g_0}} (K d\, g + K_0 d g_0) =  \int_M ( |\nabla u|_{g_0}^2 + 2 K_0 u) d\,g_0.
\]
 The derivative of this functional is:
\[
\begin{array} {cll}
{1\over 2} {{d\, F}\over {d\,t}}&= &\int_M (- 2 \triangle_0 u + K_0) {{\partial u}\over {\partial t}} d\, g_0 \\
& = & {1 \over 2} \int_M (- 2 \triangle_0 u + K_0) \triangle_g {{\partial \varphi}\over {\partial t}}  d\, g_0\\
&  = &  {1 \over 2}  \int_M  \triangle_g K \;{{\partial
\varphi}\over {\partial t}} d\, g.
\end{array}
\]
Under the Calabi flow, we have
\[
\begin{array}{lcl}
  {{d\, F}\over {d\,t}}& =&   \int_M K  \triangle_g K d\,g\\
  & = & - {1\over 2}\int_M  |\nabla K|_g^2 d\,g  \\
  & = & -{1\over 2} \int_M  |\nabla K|_{g_0}^2 d\,{g_0}\leq  0.\end{array}
\]
The last equality holds since $g$ and $g_0$ are pointwise
conformal to each other.   A well known fact is that the
Liouville functional always has a lower bound  (cf.
\cite{OPS88}). Thus, if the Calabi flow exists for time $T$ ($T$
might be equal to $\infty$), then there exists a constant $C$
independent of time $T$ such that
\begin{equation}
   \int_0^T \int_M |\nabla K|_g^2 d\, g d\,t = \int_0^T \int_M |\nabla K|_{g_0}^2 d\, g_0 d\,t < C. \label{eq: louvill bound}
\end{equation}
\end{enumerate}

  \section{Long time existence of the Calabi flow}
  In this section, we mainly use the weak compactness theorems we derived in \cite{chen942} \cite{chen981}
  and  \cite{chen970} to show that
the Calabi flow  exists for all the time. For the convenience to
the reader, we re-state the two theorems here:

\begin{theo} \cite{chen942} \cite{chen981} \cite{chen970} Let  $\{g_n, n\in {\bf N}\}$
be a sequence of conformal metrics in $\Omega$ with
 finite energy and area.  Then
there exists a subsequence of $g_n, $ a limit metric $\underline{g}$
and a finite set of points $\{p_1,p_2,\cdots p_m\}$ such that
$ g_n                   \rightharpoonup  \underline{g}\; {\rm in}\;  {H}^{2,2}_{loc}(\Omega\setminus \{p_1,p_2,\cdots p_m\})$\footnote{If $\underline{g}\equiv 0,$
the above weak converges means $g_n \rightarrow 0$ locally everywhere
except the bubble point.}. Moreover, $E(p_i) \cdot A(p_i) \geq 4\pi^2 $
where  $E(p_i)$ and $A(p_i)$ represent the amount of  energy and area
concentrated at point $p_i$ respectively. If $\underline{g} \neq 0,$
 then $\underline{g}$ has a  weak cusp singularity at each point $p_i$\footnote{Around each singular point, taking average of
 the conformal parameter over  each concentric circle, the limit
angle (if well defined) of the resulting rotationally symmetric metric is the
so called ``weak singular angle'' of the original metric.
 If the weak singular angle is $0,$ it is called a weak cusp singular point.
 The definition of weak angle could  be naturally generalized in the most natural way to the case
when the limit metric vanishes.}; the total energy
concentration could be improved as
 $ E(p_i) \cdot A(p_i) \geq 16 \pi^2;$ and
 the last  inequality is sharp.
\end{theo}

\begin{theo} (Continued from Theorem 1)\cite{chen942} \cite{chen981} \cite{chen970}.
At each bubble point $p_i,$ there exists a local re-normalization
of the metrics $h_n = \pi_n^* g_n$ such that for this new sequence
of metrics near $p_i, $
there exists  a  subsequence $\{h_{n_j},j \in {\bf N}\}$
 of $\{h_n\},$   a
finite number of bubble points $\{q_1,q_2,\cdots,q_l\}
( 0 \leq l \leq \sqrt{\frac{A(p_i)\cdot K(p_i)}{4\pi^2}})$
 with respect to the subsequence of metrics $\{ h_{n_j}\},\;$
a metric $\underline{h}$ in $  S^{2}\setminus
\{\infty, q_1,q_2,\cdots,q_l\}) $ such that: $ h_{n_j} \rightharpoonup \underline{h}$  in ${H}^{2,2}_{loc} ( S^{2}\setminus \{\infty, q_1,q_2,\cdots,q_l\}).\; $
If $\underline{h}  \equiv 0 $ (vanishing case),  then $l \geq 2 $
and $ z = 0 $ is a bubble point of $\underline{h}.\;$
If $\underline{g}$ has a non-negative weak singular angle at
$p_i,$ then a) $\underline{g}$ has a weak cusp singularity at $p_i;$
 b)  $\underline{h} \neq 0 $ and $\underline{h}$ has only weak cusp singularities (including the singular point at $z=\infty$).
\end{theo}
Suppose that $T \in (0,\infty) $ is the maximum time the flow
exists. We want to show that there is a regular metric at time
$T$ and the flow can be extended beyond time $T.\;$ For any point
$p$ in $M,$ let $\eta(x)$ be any cut off function on a unit ball
centered at $p$ w.r.t. metric $g_0.\;$ In other words,
\[
\eta(x) = \left\{ \begin{array}{l l} 1 & {\rm if}\; dist_{g_0} (x, p) < {1\over 2},\\
                                0 & {\rm if}\; dist_{g_0} (x, p) > 1,\\
                                \in (0,1) & {\rm otherwise.}
\end{array}
\right.\]
For any $\epsilon > 0, $ define an $\epsilon $ cut off function as $\eta_{\epsilon}(x) = \eta(\epsilon x)\;$
which is supported in an  $\epsilon$ ball $ B_{\epsilon}(p).\;$ Thus
\[
  \int_{B_{\epsilon}(p)}\;|\nabla \eta_{\epsilon}|^2_{g_0} \; d\, g_0 = \int_{B_{1}(p)} \;|\nabla \eta|^2_{g_0} \; d\, g_0. \]
Define a local area functional as
\[
   A_{\epsilon}(p) = \int_M \eta_{\epsilon} d\, g  = \int_{B_{\epsilon}(p)} \eta_{\epsilon} d\,g.
\]
By definition, we have
\[
  \int_{B_{{\epsilon\over 2} } } d\, g \leq A_{\epsilon} \leq \int_{B_{{\epsilon} } } d\, g. \]
 For any fixed $\epsilon, $ we have
\[\begin{array}{cll} \vert {{d A_{\epsilon} }\over {d\, t}}\vert & = & \vert \int_{B_{{\epsilon} } }
 \eta_{\epsilon} \triangle_0 K \; d\,g_0 \vert =  \vert \int_{B_{{\epsilon} } } \left(\nabla \eta_{\epsilon} \cdot \nabla K\right)_{g_0} \; d\, g_0 \vert \\ & \leq &   \left(\int_{B_{{\epsilon} } } |\nabla \eta_{\epsilon}|^2 d\,g_0\right)^{1\over 2} \cdot
\left(\int_{B_{{\epsilon} } } |\nabla K|^2_{g_0}
d\,g_0\right)^{1\over 2} \\ & \leq & C \left(\int_{B_{{\epsilon}
} } |\nabla K|^2_{g_0} d\,g_0\right)^{1\over 2}. \end{array}
\]
 Thus for any  time $t_1 <T, $ we have
\[
\begin{array}{cll}  \vert \int_{t_1}^{T} {{d A_{\epsilon} }\over {d\, t}}\; d\,t \vert
 & \leq & C \int_{t_1}^T \left(\int_{B_{{\epsilon} } } |\nabla K|^2_{g_0} d\,g_0\right)^{1\over 2}\\
&\leq & C \sqrt{T-t_1} \left(\int_{t_1}^T \int_{B_{\epsilon}(p)}
|\nabla K|^2_{g_0} d\,g\right)^{1\over 2} \\
& \leq & C \sqrt{T-t_1}.
\end{array}\]
The last inequality holds because of the inequality (\ref{eq:
louvill bound}). In other words, for any $\epsilon > 0, $ we have
(for any time $t_1 < t_2 < T$):
\[ A_\epsilon(p) \vert_{t_2} -  A_\epsilon(p) \vert_{t_1} \leq  C_1  \sqrt{t_2-t_1}
\]
or
\[
  \int_{B_{{\epsilon \over 2}} } \; d\, g(t_2) \leq  \int_{B_{\epsilon}} d\, g(t_1) +  C\sqrt{t_2 - t_1}
\]
where $C$ is a constant independent of time $t.\;$

 Now $\{ g(t)\vert 0<t< T\}$ is a 1 -parameter family of conformal metrics with finite Calabi energy
and area. Suppose that the compactness fails at least for a
subsequence $t_i \rightarrow T.\;$ According to the weak
compactness Theorem 4.1, there must exists a finite number of
points where the area function has a positive concentration.
Suppose $p$ is such a point, and $A(p) $ is the positive area
concentration. Then
\[
 A(p) = \displaystyle \overline{\lim}_{r\rightarrow 0} \displaystyle \overline{\lim}_{t_i\rightarrow T}
\int_{B_r(p)} \, d\,g(t_i) > 0.
\]
On the other hand, choose any $t <T$ and fix it for the time being:
\[
\begin{array}{lcl}
\displaystyle \overline{\lim}_{t_i\rightarrow T} \int_{B_r(p)} \,
d\,g(t_i) & \leq &  \displaystyle \overline{\lim}_{t_i\rightarrow
T} ( \int_{B_{2r}(p)} \, d\,g(t) + C \sqrt{t_i - t})\\
& \leq & \int_{B_{2r}(p)} \, d\,g(t) + C \sqrt{T - t}.
\end{array}
\]
Thus
\[
 \begin{array}{cll} A(p) & = & \displaystyle \overline{\lim}_{r\rightarrow 0} \displaystyle \overline{\lim}_{t_i\rightarrow T}
\int_{B_r(p)} \, d\,g(t_i) \int_{B_r(p)} \, d\,g(t_i)\\
 & < & \displaystyle \overline{\lim}_{r\rightarrow 0}\left(\int_{B_{2r}(p)} \, d\,g(t_0) + C \sqrt{T - t}\right)
 \\ & \leq & C\sqrt{T-t}. \end{array}
\]
Now let $t\rightarrow T, $ we have $A(p) \rightarrow 0, $  a contradiction! Thus for any sequence $t_i \rightarrow T,$ the sequence of metrics $g(t_i)$ converges to a limit metric $g(T).\;$ In other words, the conformal parameters
 $ u(t) $ of metric $g(t)$ remain uniformly
bounded from above and  below. It is then not difficult to show
that the flow actually converges to a smooth metric $g(T)$ as
$t\rightarrow T.\;$ The Calabi flow  can be extended further
beyond time $t = T.\;$ Therefore, the initial assumption the
Calabi flow exists only for a finite time is wrong and the flow
actually exists for all the time. We then have
\begin{theo} For any smooth initial metric $g_0,$ the Calabi flow exists for
all the time.
\end{theo}

\begin{rem} If one can prove that the Calabi energy is preserved along the Ricci flow (without using
the maximum principle), then  the same idea of using integral
estimates of curvature to get necessary control of the curvature
may be applied to the Ricci flow in Riemann surface as well.
\end{rem}

\section{Convergence of flow for some sequence $t_i \rightarrow \infty$}

\begin{prop}Let $\{g_i\}$ be a sequence of conformal metrics with uniformly
bounded area,  the Calabi energy and the Liouville energy. Suppose
further that one of the following holds:
\item  \begin{eqnarray} \displaystyle \lim_{i \rightarrow \infty} \int_{M} (K_{g_i} - \underline{K})^2 d\,g_i & =
& 0 ,\label{eq:conditiona}\\
 \displaystyle \lim_{i \rightarrow \infty} \int_{M} |\nabla
K_{g_i}|^2_{g_i}\; d\,g_i & = & 0.\label{eq:conditionb}
\end{eqnarray}
Then, there exists a subsequence of $g_i,$ a corresponding
sequence of  conformal transformations $\pi_i$ and a constant
scalar curvature metric $g_{\infty}$ such that $\pi_i^* g_i $
converge to $g_{\infty}$ in $H^{2,2}(M,g_0)$ in terms of the
conformal factors $\tilde{u}_i$ (If $\chi(M) < 0, $ then $\pi_i$
is trivial.), where $\pi_i^* g_i = e^{2 \tilde{u}_i} g_0.\;$
\end{prop}

We will first give a proof based on the weak compactness Theorems
4.1 and 4.2. In the case of torus and surfaces of higher genus, we
give a second proof which is more traditional. I prefer the first
one because it is more likely to be generalized in higher
dimensional manifolds, although the second proof is shorter and
cleaner.

\begin{proof} We prove this proposition by using the weak compactness
Theorems 4.1 and 4.2.
  If the compactness
fails, then there exists a finite number of bubble points $p_1,
p_2, \cdots, p_m$ such that $g_i$ has a non-trivial area
concentration in each bubble point (we follow notations in
Theorem 1 above). Consider two cases: $M = S^2$ or $M$ is a torus
or higher genus.  \\

In the first case,  $M =S^2.\;$ We can re-normalize the sequence
so that $\underline{g} \neq 0 $ (by conformal transformations).
According to Theorems 1 and 2, each bubble metric  has a weak
cusp singularity at its singular points. On the other hand,
suppose that $p$ is such a bubble point ($p$ might be any of
$p_1,p_2,\cdots, p_m$), and $D$ is a small
 disk centered at $p$ (for convenience, we assume that it is a disk
of radius 1). Now suppose the sequence of metrics can be
re-written as: $g_i = e^{2 u_i} |d\,z|^2.\;$ Here $z$ is the
coordinate variable for disk $D.\;$ According to our assumption
that there is a positive area concentration at $z=0$ (or point
$p$), thus $\displaystyle \max_{x\in D}\; u_i(x) \rightarrow
\infty.\;$ Without loss of generality, we may assume that $u_i(0)
=\displaystyle \max_{x\in D}\; \; u_i(x) \rightarrow \infty.\;$
Now re-normalize this sequence of metrics by:
\begin{equation}
   \tilde{u_i}(x) = u_i(x \;\epsilon_i) + \ln \epsilon_i,
   \label{eq:renormalize}
\end{equation}
where
\[
 \epsilon_i
   = e^{ - u_i(0)}\rightarrow 0, \qquad \forall \;| x | < {1\over {\epsilon_i}}.
\]
Note that the metric $\tilde{g}_i =  e^{2 \tilde{u}_i} |d\,z|^2 $
is just a re-normalization of $g_i,\;$ thus the scalar curvature
and area is not changed! In other words, if we denote the scalar
curvature of $\tilde{g}_i$ by $\tilde{K}_i,\;$ then
$\tilde{K}_i(x) = K_{g_i}(\epsilon_i x).\;$ Thus, for any fixed $R
> 1$ and $i$ large enough we have
\[
\begin{array}{lcl}
   \int_{|x| < R}\; {{(\triangle \tilde{u}_i)^2}\over{e^{ 2\tilde{u}_i}}} \;|d\,z|^2
  & = & \int_{|x|< R} \tilde{K}_i^2 \;d\,\tilde{g}_i
  \\
  & = & \int_{|x| < \epsilon_i R} K_{g_i}^2 \;d\;g_i \\
  & \leq & \int_{D} \; K_{g_i}^2 \;d\,g_i <C_1,
  \end{array}
\]
and
\[
\begin{array}{lcl}
\int_{|x| < R}\; e^{ 2\tilde{u}_i} \;|d\,z|^2 & = &   \int_{|x| <
\epsilon_i R}\;e^{ 2 u_i} |d\,z|^2 \\
& \leq & \int_D\; d\;g_i \leq  C_2.
\end{array}
\]
By definition, we have $\tilde{u_i}(x)\leq 0$ for any $|x| < R.\;$
Therefore, it is not difficult to choose  a subsequence of
$\tilde{g}_{i}$ which converges in every fixed disk $|x|< R.\;$
Suppose the limit metric is $\tilde{g}.\;$ Then $\tilde{g}$ is a
constant scalar curvature metric (Equation (\ref{eq:conditiona})
or (\ref{eq:conditionb}) implies this) in the Euclidean plane
with finite area. Thus, $\tilde{g}$ is a smooth metric in $S^2$
with positive constant curvature. This contradicts the earlier
assertion that any bubble metric must only have weak cusp
singularities (since the scalar curvature must be negative near a
cusp singular point). Therefore, there is no bubble point after a
possible re-normalization of conformal transformation. In other
words, there is a subsequence of $g_i$ which converges weakly in
$H^{2,2}(S^2)$ (up to conformal transformation group) to a metric
$g_{\infty}.\;$ Equation (\ref{eq:conditiona}) or
(\ref{eq:conditionb}) implies that $g_{\infty}$
has constant scalar curvature.\\

In the second case, $M$ is either a  torus or  surface of higher
genus. We can not re-normalize like in $S^2$ case.
 There are two cases
to handle,  the  case when $\underline{g} \neq 0,\;$ and the case
when $\underline{g} \equiv 0.\;$  In the first case when
$\underline{g} \neq 0,\;$ one can argue like in $S^2:\; $ the
bubble metric if existed, must be a round sphere on one hand; and
must have one cusp singularity on the other hand. This is
impossible since the scalar curvature near cusp singular point
must be negative. Therefore there are no bubble point in the first
case and the limiting metric must have constant scalar curvature
metric (cf. the assumption  (\ref{eq:conditiona}) or
(\ref{eq:conditionb})).  In the second case, we have
$\underline{g} \equiv 0.\;$ We will show that this is impossible
by drawing a contradiction. First, it
 is a ghost vertex \footnote{We can apply Theorem 4.2 iteratively to any bubble point arisen
  from taking weak limit of some subsequence
 of $\{g_i\}$, we eventually obtain  a
  "bubbles on bubbles" phenomena. The limit
 of metrics at each stage of "blow-up" is regarded as a "vertex" in the limit tree structure of $\{g_i\}.\;$
 If the limit metric vanishes
 identically, it is then called a "ghost vertex." It was proved in \cite{chen942}  that there are only a finite
 number of ghost vertexes. } in the tree decomposition of the limit of
$\{g_i\}.\;$ Theorems 4.1 and 4.2 imply that the limit of
$\{g_i\}$ can be decomposed in a tree structure where each vertex
represents a limit of  a subsequence of $g_i$ under a local
re-normalization (cf. equation (\ref{eq:renormalize})).  The
assumption (\ref{eq:conditiona}) or (\ref{eq:conditionb}) ensures
that each vertex (non-ghost) corresponds to a metric with
constant scalar curvature and finite area. Therefore,  each
non-ghost vertex must be a sphere! Consider the total Euler
character in the limit tree structure, it should be non-positive
since this is torus or high genus surface. On the other hand, the
contribution from each non vanishing vertex is always positive
(since they are $S^2$); thus the total Euler character in the
remaining ghost vertexes must be strictly negative (actually less
than $- 4\pi$ if there is one $S^2$ in the limit!).  Moreover,
the total area concentration is $0$ at each ghost vertex.
Therefore, the total energy concentration in the ghost vertex
must be infinite by the Schwartz inequality:
\[
  (- 4\pi)^2  \leq (\int_{\Sigma} K_{g_i})^2 \;d\,g_i  \leq  \int_{\Sigma} 1 \cdot \int_{\Sigma} K_{g_i}^2\;d\,g_i,
\]
where $\Sigma$ denotes the total collection of "ghost vertexes."
$ \int_{\Sigma} K_{g_i}^2\;d\,g_i \rightarrow \infty $ because $
\int_{\Sigma} 1 \rightarrow 0.\;$ However, the total energy of
the tree structure is finite, this is a contradiction! Thus,
$g_i$ must converge to some constant scalar curvature metric in
both torus and surfaces of high genus.
\end{proof}

Now we give the second proof on the case of torus or surface with
high genus.
\begin{proof} Write $g_i = e^{2 u_i} g_0.\;$ Note that the Calabi
energy, Liouville energy and area are uniformly bounded along
this sequence of metrics. Then we have
\begin{equation}
\displaystyle \int_M {{(-\triangle_0 u_i + K_0)^2}\over {e^{2
u_i}}} \;d\,g_0  \leq  C, \label{eq:calabienergy}
\end{equation}
and \begin{equation} - C \leq  \displaystyle \int_M \mid \nabla
u_i\mid^2_{g_0} \;d\,g_0 + \displaystyle \int_M 2 K_0 u_i\;
d\;g_0  \leq  C. \label{eq:liouvilleenergy}
\end{equation}
Here $C$ is some uniform constant in this proof, and its value
may change from line to line. Since area is fixed along the flow:
$ \displaystyle \int_M e^{2 u_i}\;d\,g_0  = \displaystyle \int_M
\;d\,g_0. $ It follows that
\[
  \displaystyle \int_M u_i d\,g_0 \leq C.
\]
Since the underlying surface is torus or surface of higher genus,
we may assume (without loss of generality) that $K_0$ is a
non-positive constant. Combining the previous inequality and the
equation (\ref{eq:liouvilleenergy}), we obtain
\[
  \displaystyle \int_M \mid \nabla u_i \mid_{g_0}^2 \;d\,g_0 \leq C.
\]
Following from the Moser-Trudinger Inequality, we have
\[
  \displaystyle \int_M\; e^{3 u_i} \;d\,g_0 \leq C.
\]
Combining this with equation (\ref{eq:calabienergy}), we arrive at
\[
  \displaystyle \int_M \;(-\triangle_0 u_i + K_0)^{6\over 5} \;d\,g_0 \leq
  C.
\]
Note that $ \displaystyle \int_M e^{2 u_i}\;d\,g_0  =
\displaystyle \int_M \;d\,g_0. $ It then follows that
\[
  \| u_i \|_{H^{2, {6\over 5}} (M,g_0)} \leq  C.
\]
In particular, $\| u_i\|_{L^\infty} \leq C.\;$  In view of
Theorems 4.1 and 4.2, there is no concentration point for this
sequences. It follows that $\{g_i\}$ always converges by sequence.
The condition (\ref{eq:conditiona}) or (\ref{eq:conditionb})
implies that the limit metric has constant scalar curvature
metric.
\end{proof}

\begin{defi}
For any metric $g$, define $\lambda_1(g)$ to be the first
eigenvalue  of the complex Laplacian
 operator with respect to the metric $g.\;$ This is a well defined
map  from the space of metrics to the positive real line.
\end{defi}

Clearly, $\lambda_1(g)$  is invariant under conformal
transformation (where the metric is viewed in a different
coordinate, but not re-scaled!).

\begin{prop} Let $\{g_i\}$ be a sequence of conformal metrics with fixed area
(normalized so that total area is $ 2 \pi |\chi(M)|$ except in a
torus where one could normalized area to be any size). Suppose
that $g_i = e^{2 u_i} g_0$ weakly converge to a constant scalar
curvature metric $g_{\infty}$ in $H^{2,2}(M,g_0)$ in terms of
conformal parameter $u_i.\;$ Then $\displaystyle
\lim_{i\rightarrow \infty} \lambda_1(g_i) = \lambda_1 > 0
$\footnote{$\lambda_1(g_{\infty}) = 1$ if $M = S^2.\;$ }. As an
immediate corollary, if $\int_M |\nabla K_{g_i} |^2_{g_i}\; d\,
g_{i} \rightarrow 0,\;$ then $\int_M |K_{g_i} - \underline{K}|^2
d\, g_{i} \rightarrow 0.\;$

\end{prop}
\begin{proof}
 Suppose $g_{\infty} = e^{ 2 u_{\infty}} g_0.\;$ Then $u_i
\rightarrow u_{\infty}$ in $H^{2,2}(M,g_0).\;$ Thus $u_i
\rightarrow u_{\infty}$ in $C^{\alpha}$ for any $0 < \alpha <
1.\;$ In other words, $ e^{ 2 u_i} \rightarrow e^{ 2 u_{\infty}}
$ uniformly. The later in turn implies that first eigenvalues of
$g_{i}$ must converge to that of $g_{\infty}.\;$ The last
statement of the proposition just follows from the Poincare
inequality.
\end{proof}

\begin{rem} The first eigenvalue of the evolved metrics converges
to a positive constant implies that the Sobolev constant for the
evolved metric is uniformly bounded from below.
\end{rem}
\begin{theo}Let $\{g(t) \vert 0\leq t < \infty\}$ be a one-parameter family of metrics under the Calabi flow.
Then for any sequence of numbers $t_i \rightarrow \infty,$ there
exists a subsequence (denoted again by $\{t_i\}$) such that
$g(t_i)$ converges weakly to a constant scalar curvature metric
up to conformal transformation.
\end{theo}
\begin{proof} Since the Calabi flow exists for long time,  the
inequality (\ref{eq: louvill bound}) holds for time $T = \infty :$
\[
   \int_0^{\infty} \int_M |\nabla K|_g^2 d\, g d\,t = \int_0^{\infty} \int_M |\nabla K|_{g_0}^2 d\, g_0 d\,t < C
\]
for some constant $C > 0.\;$ Therefore there exists at least a
sequence of numbers $t_i \rightarrow \infty$ such that
\[
  \displaystyle \lim_{t_i \rightarrow \infty} \int_M |\nabla K_{g(t_i)}|_{g(t_i)}^2 d\, g(t_i) d\,t =    \displaystyle \lim_{t_i \rightarrow \infty} \int_M |\nabla K_{g(t_i)}|_{g_0}^2 d\, g_0 = 0.
\]
Proposition 5.1 implies that there exists a subsequence (denoted
again by $t_i$), a subsequence of conformal transformations
$\pi_i$ and a constant scalar curvature metric $g_{\infty}$ such
that $\pi_i^* g(t_i) \rightarrow g_{\infty}$ weakly in
$H^{2,2}(M,g_0).\;$ Now Proposition 5.3 implies that
\[\begin{array}{lcl}
\displaystyle \lim_{t_i \rightarrow \infty} \lambda_1(g(t_i)) & =
& \displaystyle \lim_{t_i \rightarrow \infty} \lambda_1(\pi_i^*
g(t_i) ) \\
& = & \lambda_1(g_{\infty}) > 0,\end{array} \]
 and
\[
\begin{array}{lcl}
\displaystyle \lim_{t_i \rightarrow \infty} Ca(g(t_i)) & = &
\displaystyle \lim_{t_i \rightarrow \infty} Ca(\pi_i^* g(t_i)) \\
& = & \displaystyle \lim_{t_i \rightarrow \infty} \int_M
|K_{g(t_i)} - \underline{K}|^2 d\, g(t_{i}) = 0. \\
\end{array}
\]

However, the Calabi energy decreases monotonely under the Calabi
flow. The above inequality implies that
\[
\displaystyle \lim_{t \rightarrow \infty} Ca(g(t))
= \displaystyle \lim_{t \rightarrow \infty} \int_M |K_{g(t)} - \underline{K}|^2 d\, g(t) = 0.
\]
For any sequence $t_i \rightarrow \infty, $ Proposition 5.1 again
implies that there exists a subsequence (denoted again by
$t_i$),  a sequence of conformal transformation $\pi_i, $ and a
constant scalar curvature metric $g_{\infty}$ such that $\pi_i^*
g(t_i) \rightarrow g_{\infty}$ weakly in $H^{2,2}(M,g_0). \;$
\end{proof}
 Notice that
the constant scalar curvature metric $g_{\infty}$ so obtained
depends on the sequence chosen. However, all constant scalar
curvature metrics in $M^2$ with same area have the same spectrum,
in particular, the same first and second eigenvalues. Therefore,
if one denotes $\lambda_1(g_\infty), \;\lambda_2(g_\infty)$ as the
first and second eigenvalues of the limit constant curvature
metric, it is then well defined (independent of a sequence
chosen). As a matter of fact, on $S^2$ we have
$\lambda_1(g_\infty) = 1 $ and $\lambda_2(g_\infty)=3.\;$ This
leads to the following corollary:
\begin{cor} Let $\{g(t) \vert 0\leq t < \infty\}$ be a one-parameter family of metrics under the
 Calabi flow. For any $\epsilon > 0,\;$   for $t$ large enough,  the eigenvalues of $g(t),$  either
 are between $(1-\epsilon,1+\epsilon)$ or are bigger than $2.\;$ Let $\Lambda_{\rm first}$
 denote the eigenspace of $g(t)$ corresponds to the eigenvalues between  $1-\epsilon$ and $1+\epsilon.\; $ Then
 $\Lambda_{\rm first}$  converges to the first eigenspace of some constant scalar curvature metric in the limit.

\end{cor}

\section{Uniqueness of the limit at $t \rightarrow \infty$
for different sequences }

\begin{prop} There exists a small positive number $\alpha $ and a big constant $C$ such that
\[ \int_M (K - \underline{K})^2 d\,g(t) < C \; e^{- \alpha \; t}, \qquad \forall\; t > 0.\]
\end{prop}
\begin{proof} Let $\epsilon(t)^2 = \int_M (K - \underline{K})^2
d\,g(t)$ and let $C$ denote any generic constant. Then
$\epsilon(t) \rightarrow 0$ as $t\rightarrow \infty.\;$ We want
to show that $\epsilon(t)$ decays exponentially fast (note that
$K^{,zz} = { 1 \over F^2} K_{,zz}$ in the following calculation):

\begin{eqnarray} & & {d\over {d\,t}} \int_M (K - \underline{K})^2 d\, g(t) \nonumber \\
 & = & -\int_{M} K_{,zz} K^{,zz} d\,g(t) \nonumber\\
  & = & - \int_M (\triangle_g K)^2 d\,g(t) + \int_M {K\over 2} |\nabla K|^2_g d\,g \nonumber \\
  & = & - \int_M (\triangle_g K)^2 d\,g(t) + \int_M {{K - \underline{K}}\over 2} |\nabla K|^2_g d\,g + {{\underline{K}}\over 2} \int_M  |\nabla K|^2_g d\,g
   \nonumber \\
  & \leq & - \int_M (\triangle_g K)^2 d\,g(t) + {1\over 2} (\int_M (K - \underline{K})^2 d\,g)^{{1\over 2}}\cdot (\int_M |\nabla K|^4_g d\,g)^{{1\over 2}}
   + {{\underline{K}} \over 2} \int_M |\nabla K|^2_g\; d\,g  \nonumber\\
  & = & - \int_M (\triangle_g K)^2 d\,g(t) + C \epsilon(t)^{{1\over 2}}\times (\int_M |\nabla K|^4_g d\,g)^{{1\over 2}}
   + {{\underline{K}} \over 2} \int_M |\nabla K|^2_g d\,g. \label{eq:expodecay}
\end{eqnarray}

Now we need to consider two cases: $\underline{K} \leq 0$ or
$\underline{K} > 0.\;$ The first case is easier, while the second
case is more delicate. \\

Consider the first case.  If $\underline{K} \leq 0, $ then there
exists some positive constant $\alpha$ and $C$ such that
\[
\begin{array}{lcl} & &  {d\over{d\,t}} \int_M (K - \underline{K})^2 d\,
g(t)\\  & \leq & - \int_M (\triangle_g K)^2 d\,g(t) + C
\epsilon(t)\times (\int_M |\triangle_g K|^2_g d\,g)
-  {|{\underline{K}|} \over 2}  \displaystyle \int_M  |\nabla K|^2_g d\,g \\
& \leq & - C \int_M (\triangle_g K)^2 d\,g(t) \leq -\alpha \int_M (K -\underline{K})^2 d\,g(t)
\end{array}
\]
for some positive constant $\alpha > 0.\;$ Thus, we have
\[
 \int_M (K - \underline{K})^2 d\,g \leq  C e^{-  \alpha t}, \qquad {\rm for\; any \;} t > 0.
\]

Next we consider the second case.  If $\underline{K} > 0,$ then
our normalization yields that $\underline{K} =1.\;$ We can
rewrite equation (\ref{eq:expodecay}):
\begin{eqnarray} & &
{d\over {d\,t}} \int_M (K - \underline{K})^2 d\, g(t) \nonumber \\
& \leq & - ( 1 - C \epsilon(t))  \int_M (\triangle_g K)^2 d\,g(t)
   - {{1} \over 2} \int_M (K - \underline{K}) \triangle_{g(t)} K\; d\,g \nonumber
   \\ & = & -(1 - C \epsilon(t))  \int_M (\triangle_g K)^2 d\,g(t) + {1\over 2}  \int_M (\triangle_g K)^2 d\,g(t)
   +  {1\over 2}  \int_M (K - \underline{K})^2 d\,g(t) \nonumber
   \\
     & = & -({1\over 2} - C \epsilon(t))  \int_M (\triangle_g K)^2 d\,g(t) +
      {1\over 2}  \int_M ( K - \underline{K})^2 d\,g(t)
      \label{eq:expodecay1}.
\end{eqnarray}
Since the first eigenvalue of $g(t)$ converges to $1,$ the above
inequality seems to give us little control on the decay rate of $
{d\over {d\,t}} \int_M (K - \underline{K})^2 d\, g(t).\; $
However, by using Kazdan-Warner condition, we can still prove the
above estimate in $S^2.\;$ Notice that Kazdan-Warner condition
implies
\[
 \displaystyle \int_M X(K- \underline{K}) d\;g(t)= \displaystyle \int_M (K- \underline{K}) \triangle_{g(t)}
 \theta_X d\;g(t) = 0,
 \]
 where $X$ is any holomorphic vector field in $S^2,\;$ and $\theta_X $ is the potential of the vector
 field $X$ with respect to the metric $g(t)$\footnote{For any holomorphic vector field $X$ on $S^2,$ and for any metric $g,$ one can
 define the  potential function $\theta_X$  (up to addition of some constants) as the following:
 \[
   \displaystyle \int_M \;X(\psi)\; d\,g = -  \displaystyle \int_M
   \;\theta_X \triangle_g \psi\; d\,g,
 \]
 where $\psi$ is any smooth test function.  If $g $ is  a
 constant scalar curvature metric, then $\theta_X$ is an eigenfunction  of $g$ with eigenvalue
 $1,\;$i.e.,
 \[
 \triangle_{g} \theta_X = - \theta_X.
 \]
Recall $\eta(S^2)$ denote the lie algebra of all holomorphic
vector fields on $S^2.\;$ Then
\[ \{\triangle_{g} \theta_X\mid \forall\; X \in \eta(S^2)\} = \{ \theta_X
\mid \forall\; X \in \eta(S^2)\}
\]
is the first eigenspace of $g.\;$ }.
  The Kazdan-Warner
condition can be re-written as
\begin{equation}
\displaystyle \int_M\;(K- \underline{K}) \triangle_{g(t)} \theta_X
e^{2u} d\,g_{\infty} =0. \label{eq:Kazdan-warner}
\end{equation}


   Note that $g(t)$ converges to some constant scalar
curvature metric in $S^2$ in the $H^{2,2}$ sense. This should be
enough for our purpose.  For $t$ large enough, we can choose a
constant scalar curvature metric $ g_{\infty}$ (which may depends
on $t$) such that
\[
g(t) = e^{2 u} g_{\infty}
\]
and
\[
  \|u\|_{H^{2,2}(S^2,  g_{\infty})} \rightarrow 0
\]
as $t \rightarrow 0.\;$  This means that the eigenvalue of $g(t)$
converges to those  of  some constant scalar curvature metric
(whose first, second eigenvalues are $1,3$ respectively). As
Corollary 1 implies, the eigenspace $\Lambda_{\rm first}$ of
$g(t)$ corresponds to eigenvalues between $1-\epsilon$ and
$1+\epsilon$ converges to the first eigenspace of $g_{\infty}.\;$

If we decompose $K - \underline{K}$ into two components:
 \[
   K - \underline{K} = \rho + \rho^{  \perp}
 \]
 where
 \[
\rho \in \Lambda_{\rm first}(g(t)),\qquad {\rm and} \qquad
\rho^{\perp} \perp \Lambda_{\rm first} (g(t)).
 \]
Then the Kazdan-Warner (\ref{eq:Kazdan-warner}) condition implies
that \footnote {Note that $\{\triangle_{g(t)} \theta_X \mid
\forall X \in \eta(S^2) \}$ and $\Lambda_{\rm first}(g(t))$ both
converge to the first eigenspace of same constant scalar
curvature metric $g_{\infty}$ in the limit as $t\rightarrow
\infty.\;$ Thus the Kazdan-Warner (\ref{eq:Kazdan-warner})
condition implies that the projection of  $K- \underline{K}$ into
$\Lambda_{\rm first}(g(t))$ is very small compared to the $L^2$
norm of $K- \underline{K}.$}:
\[
   \|\rho\|_{L^2(g(t))} \leq \epsilon\; \|K -
   \underline{K}\|_{L^2(g(t))},
 \]
 where $\epsilon \rightarrow 0$ as $t \rightarrow \infty.\;$
 Plugging this into equation (\ref{eq:expodecay1}), we have
 (note that $\underline{K} =1$ in the following calculation):
\[
\begin{array}{lcl}& & {d\over {d\,t}} \int_M (K - \underline{K})^2 d\, g(t) \\
& \leq & ({1\over 2} - C \epsilon(t))  \int_M (\triangle_g K)^2
d\,g(t) +  {1\over 2}  \int_M ( K - \underline{K})^2 d\,g(t) \\
   & \leq & - ( {1\over 2} - C \epsilon(t) ) 2^2 (1-\epsilon)
    \int_M (K - \underline{K})^2 d\,g(t) + {1 \over 2}  \int_M
(K-\underline{K})^2 d\,g(t)
   \\
   & \leq & -\alpha  \int_M (K-\underline{K})^2 d\,g(t),
\end{array}
\]
where $\alpha =  ( {1\over 2} - C \epsilon(t) ) 4 (1-\epsilon)- {1
\over 2}.\; $ Choose $\epsilon$ small enough and $t$ large enough
(note that $\displaystyle \lim_{t\rightarrow \infty} \epsilon(t)
= 0$), we have $\alpha
> 0.\;$ It follows that we can easily see  the exponential
decay of $L^2$ norm of $K - \underline{K}$:
\[
\int_M (K - \underline{K})^2 d\,g \leq  C e^{-  \alpha t}, \qquad
{\rm for\; any \;} t > 0.
\]
\end{proof}
 The key
estimate we used in this calculation is that the spectrum of the
evolved metrics $g(t)$ converges to that of a constant scalar
curvature metric. All of the constant scalar curvature metrics in
$S^2$ (with same area) are isometric to each other. Thus any
geometric norm, such as $L^2$ norm of $\nabla K$ or any higher
order derivatives, can be computed and proved to converge
exponentially in a similar fashion,
 just as Chrusciel did in his original paper. Here we want to prove the uniqueness directly
 using this idea in infinite dimensional space.  The point of view one should
adopt is that the one parameter family of metrics $\{g(t)| 0 \leq
t <\infty\}$ represents a curve in $\cal {H}.\;$ If this curve
runs to $\infty$ in $\cal H$ in terms of geodesic distance, then
Calabi flow diverges. If this curve stays within a bounded domain
in $\cal H,$
 then this corresponds to the case that each subsequence converges , but converges to different
limit as $t \rightarrow \infty.\;$ However, the case we are in is the best possible: the Calabi
flow represents a "Cauchy" sequence under this Riemannian metric. Thus, the limit metric must be
unique. To verify this, for any $t_i >  s_i \rightarrow \infty,$ we want to show the length
of curves $\{g(t) \vert s_i < t < t_i\}$ tends to $0$  fast.  Denote the distance of
two metrics $g(t_i)$ and $g(s_i)$ by $ d(g(t_i), g(s_i)).\;$Then
\[
\begin{array}{cll}
d(g(t_i), g(s_i)) & \leq & \int_{s_i}^{t_i} \sqrt{\int_M ({{\partial \varphi}\over {\partial t}})^2 \; d\; g(t)} d\,t
\\ & = & \int_{s_i}^{t_i} \sqrt{\int_M (K -\underline{K})^2 d\,g(t)} \; d\,t \\
& \leq & \int_{s_i}^{t_i} C e^{-\alpha t} \; d\,t\\
 & \leq &  C {{e^{-\alpha s_i} - e^{-\alpha  t_i}}\over {\alpha}} \rightarrow 0.
\end{array}
\]
Thus the limit as $t \rightarrow \infty$ must be unique. We then
proves the following
\begin{theo} There exists a unique constant scalar curvature metric $g_{\infty} \in \cal H$ such
that the Calabi flow converges to this metric exponentially fast
in terms of geodesic distance.
\end{theo}
\begin{rem}  $g_{\infty}$ could still be a constant scalar curvature metric
concentrated on one point. In the following subsection, we
 will show that
this concentration do not  occur along the Calabi flow.
\end{rem}
\subsection{Convergence of conformal factors $e^{2 u}$ as $t \rightarrow \infty$}

\begin{theo} For any metric $g_0$ in Riemann surface, the Calabi  flow exists
for all the time. For every sequence $t_i \rightarrow \infty,$
there exists a subsequence $t_i$ such that $g(t_i) $ converge to
a constant scalar curvature metric. And this limiting  constant
scalar curvature metric is independent of sequence chosen.
\end{theo}
In light of Theorem 6.2 and the remark following it, we only need
to show that conformal parameters does not concentrate in any
point; or to show it is always bounded.\\

\begin{proof} We already know that there exists a family of conformal
transformations $\pi_i$ such that  $\pi_i^* g(t_i) $ converges to
a constant scalar curvature metric $g_{\infty}.\;$ Then the set
$\{\pi_i|i = 1,2, \cdots \}$ must be compact \footnote{This is
equivalent to say $|u_i| \leq C $ for some uniform constant $C$
for all $i \geq 1.\;$ Here $u_i$ is the conformal parameter of
$g_i,$ i.e., $g_i = e^{2 u_i} g_0.\;$ }. Otherwise, suppose  that
$\{\pi_i|i = 1,2, \cdots \}$ is a non-compact family of conformal
transformations. By a direct calculation, one should yield:
\[
   d(g_{\infty}, \{\pi_i^{-1}\}^* g_{\infty}) \rightarrow \infty.
\]
On the other hand,
\[
  0 = \displaystyle \lim_{i \rightarrow \infty} d(g_{\infty}, \pi_i^* g(t_i))
   = \displaystyle \lim_{i \rightarrow \infty} d(g(t_i), \{\pi_i^{-1}\}^* g_{\infty}).
\]
 By a triangle inequality for distance function, we have
\[
\begin{array}{lcl}
  \displaystyle \lim_{i \rightarrow \infty} d( g_{\infty}, g(t_i)) &  \geq &
   \displaystyle \lim_{i \rightarrow \infty} ( -  d(  \{\pi_i^{-1}\}^* g_{\infty}, g(t_i))
     +  d(g_{\infty}, \{\pi_i^{-1}\}^* g_{\infty})) \\
     & = & \infty.
\end{array}
\]
However, $g_{\infty}$ is invariant under the Calabi flow,
$g(t_i)$ is in the image of $g_0$ under Calabi flow at time $t =
t_i.\;$ By Theorem C,   the  distance in $\cal H$ decreases under
the Calabi flow. Consequently, we have
\[
 \infty = \displaystyle \lim_{i \rightarrow \infty} d( g_{\infty}, g(t_i))
\leq d( g_{\infty}, g_0) < \infty.
\]
This is a contradiction! Thus the set of conformal transformations
$\{\pi_i\}$ must be compact and the conformal factor must be
bounded. \end{proof}

\section{Future questions}
In this section, we list some problems we think they might be interesting, and
at the same time, accessible.
Suppose that ${\cal A }$ is a subset of the space of subset of $\cal H.\;$ Define
the diameter of this set as the maximal distance under this metric
defined by \cite{Ma87} \cite{Semmes92} and \cite{Dona96}.
\begin{q} Considering metric in a Riemann surface, suppose ${\cal A} = \{ g \vert A(g) < C, Ca(g) < C\}$ and diameter of $\cal A$ is finite, is  $\cal A$ compact?
\end{q}
In \cite{chen942}, we give a weak compactness theorem of $\cal
A.\;$ For a sequence of metrics with finite energy and area, the
weak compactness fails because large area concentrates at some
isolated points. We  ask   if this distance of infinite
dimensional space necessary approaches to $\infty $ once the
concentration of area occurs in some isolated points.

\begin{q} What is the relationship  between the
distance between any two K\"{a}hler metrics in a fixed K\"{a}hler
class (defined on \cite{Ma87}, \cite{Semmes92} and \cite{Dona96})
and the Gromov's Hausdorff distance between any two Riemannian
metrics. Certainly, the later one is more general, but one wants
to understand the difference and similarity when restricted to the
space of K\"{a}hler metrics. The question above is an attempt in
this direction.
\end{q}
\begin{q} Can one prove the lower bound of Mabuchi energy directly (without
appealing to the Mabuchi-Bando's theorem)? On a Riemann Surface,
 the Mabuchi energy takes the
form\footnote{The explicit formulation of the Mabuchi energy in
general K\"ahler class was first given by  Tian, cf. \cite{tian98}.}:
\[
  Ma(\varphi)  = \int_M \ln (1 + \triangle \varphi) (1 +\triangle \varphi) - {1\over 2} \underline{K} |\nabla \varphi|^2
      - (K_0 - \underline{K} ) \varphi
\]
for any $\varphi$ such that $ 1 + \triangle \varphi > 0 $ on $M$.
Here all of the norm are taken w.r.t. the metric $g_0.\;$ It is
easy to see that the Euler-Lagrange equation for this functional
is $\;K - \underline{K} = 0,\;$ if we recall the formula
(\ref{eq:curvature})  for scalar curvature in this setting. The
question we want to ask is if there exists a universal constant $
C > 0$ such that for any $\varphi,$ the following inequality
holds:
\begin{equation}
  \int_M \ln (1 + \triangle \varphi) (1 +\triangle \varphi) - {1\over 2} \underline{K}  |\nabla \varphi|^2 - (K_0 - \underline{K} ) \varphi > -C ?
\label{eq:Mabuchi-bound}
\end{equation}
We know this inequality holds since there always exists a
K\"{a}hler-Einstein metric in Riemannian surface and a  theorem
of Mabuchi-Bando implies  that the Mabuchi energy bounded from
below if there is a K\"{a}hler-Einstein metric. The question we
ask here is if we can get this inequality without appealing to
this theorem of Mabuchi and Bando. If we can prove this in
Riemannian surface, then there is at least some hope that we
might be able to get a lower bound  in higher dimension
in some cases without knowledge of K\"{a}hler-Einstein metrics. \\

In a way, this should be similar to Moser-Trudinger-Onofri inequality
and one should be able to get this with pure analytic tools.
\end{q}

\begin{q} The proof of lower bound of the Liouville energy makes a strong
use of Trudinger inequality and it looks quite hard to get a new
proof without assuming uniformization theorem. However, it might
be a lot easier to prove the same statement when restricted to
those conformal metrics with bounded the Calabi energy and area.
This will be sufficient for our purpose since the Calabi flow
decreases the Calabi energy and preserves the total area.
\end{q}
\begin{conj} In a topological 2-sphere  $M,$ suppose $g_i = e^{ 2 u_i} g_0$ is a family
 of conformal metrics with bounded area and energy. Suppose the following
holds
\begin{enumerate}
\item
\[ \displaystyle \lim_{i \rightarrow \infty} \int_{M} |\nabla K_{g_i}|^2 d\,g_0 = 0.
\]
\item There exists only one bubble point $p \in S^2,$ and $u_i \rightarrow  -\infty$ in any
compact subset $\Omega \subset  M \setminus \{p\}.\;$

\end{enumerate}

we conjecture that there exists a flat metric $g_{\infty} = e^{2 u_{\infty}} g_0$
in $M\setminus \{p\}$ and a subsequence of $g_i$ such that the following holds
\[
  u_i - c_i \rightharpoonup u_{\infty} \;\;{\rm in}\; H^{2,2}(\Omega)
\]
where $\Omega $ is any compact subset of $M\setminus \{p\}, $ and $c_i
= u_i(q)$ for any fixed point $q \in \Omega \subset M \setminus \{p\}.
$
\end{conj}
\begin{rem} This conjecture implies that $M$ is the standard $S^2.\;$
This conjecture, together with the lower bound of the Liouville
energy,  will re-prove the uniformization theorem via the Calabi flow.
\end{rem}

\begin{rem} There is a local version of the above conjecture which
had been worked out by Yan Yan Li \cite{Li99} and it is
interesting to compare the two problems.
\end{rem}



\noindent {\rm Added after the proof:} We notice the recent work
of M. Struwe, following a similar idea of this paper, giving a
more concise proof to the convergence of the Calabi flow from
analytic point of view; moreover,  he gives a new proof  to the
Ricci flow in Riemannian surfaces (using the idea of integral
estimate only). One key lemma he proves is that the Calabi energy
is preserved and eventually improved in the Ricci flow in Riemann
surface.

\noindent Department of Mathematics, Princeton University,
          Princeton, NJ 08544,\\
\noindent  xiu@math.princeton.edu
\end{document}